\documentclass{article}
\usepackage{spconf}
\usepackage{amsmath}
\usepackage{array}
\usepackage{bbold}
\usepackage{graphicx, graphics, theorem, times, amsfonts, amsmath, amssymb, cite}
\usepackage{tikz}
\usepackage{algorithmic}
\usepackage{algorithm}
\usetikzlibrary{shapes,arrows}
\usepackage{mathtools}
\usepackage{multirow}
\usepackage{cleveref, subcaption}
\usepackage{bm}
\usepackage{dcolumn}% Align table columns on decimal point
 \usepackage{amssymb}
\usepackage{dsfont}
\usepackage{mathdots}
\usepackage{color}
\usepackage{dcolumn}% Align table columns on decimal point

\newenvironment{myproof}[1][$\!\!$]{{\noindent\bf Proof #1: }}
                         {\hfill\QED\medskip}

\newcounter{excercise}
\newcounter{excercisepart}

\definecolor{pennblue}{cmyk}{1,0.65,0,0.30}
\definecolor{pennred}{cmyk}{0,1,0.65,0.34}
\definecolor{mygreen}{rgb}{0.10,0.50,0.10}

\def \reals    {{\mathbb R}}

\newcommand{\argmax}{\operatornamewithlimits{argmax}}
\newcommand{\argmin}{\operatornamewithlimits{argmin}}

\def\ccalB{{\ensuremath{\mathcal B}}}

\def\ccalG{{\ensuremath{\mathcal G}}}

\def\ccalI{{\ensuremath{\mathcal I}}}

\def\ccalM{{\ensuremath{\mathcal M}}}

\def\ccalP{{\ensuremath{\mathcal P}}}

\def\ccalV{{\ensuremath{\mathcal V}}}

\def\ccal0{{\ensuremath{\mathcal 0}}}
\def\bbA{{\ensuremath{\mathbf A}}}

\def\bbD{{\ensuremath{\mathbf D}}}

\def\bbH{{\ensuremath{\mathbf H}}}
\def\bbI{{\ensuremath{\mathbf I}}}
\def\bbJ{{\ensuremath{\mathbf J}}}

\def\bbL{{\ensuremath{\mathbf L}}}

\def\bbU{{\ensuremath{\mathbf U}}}
\def\bbV{{\ensuremath{\mathbf V}}}

\def\bbh{{\ensuremath{\mathbf h}}}

\def\bbn{{\ensuremath{\mathbf n}}}

\def\bbu{{\ensuremath{\mathbf u}}}
\def\bbv{{\ensuremath{\mathbf v}}}

\def\bbx{{\ensuremath{\mathbf x}}}
\def\bby{{\ensuremath{\mathbf y}}}

\def\bb0{{\ensuremath{\mathbf 0}}}

\newcommand{\QED}{\hfill\ensuremath{\blacksquare}}

\newtheorem{mytheorem}{\hspace{-1pt}\bf Theorem}
\newtheorem{mydefinition}{\hspace{-1pt}\bf Definition}

\newtheorem{mylemma}{\hspace{-1pt}\bf Lemma}
\newtheorem{myproposition}{\hspace{-1pt}\bf Proposition}
\newtheorem{remark}{\hspace{-1pt}\bf Remark}

\ninept

\newcommand{\abs}[1]{\lvert #1 \rvert}

\def\lovasz{Lov\'asz }

\title{A Digraph Fourier Transform {With Spread Frequency Components}}
\name{Rasoul Shafipour$^{\dag}$, Ali Khodabakhsh$^{\ddagger}$, Gonzalo Mateos$^{\dag}$, and Evdokia Nikolova$^{\ddagger}$}
\address{$^{\dag}$Dept. of Electrical and Computer Engineering, University of Rochester, Rochester, NY, USA\\
	$^{\ddagger}$Dept. of Electrical and Computer Engineering, University of Texas at Austin, Austin, TX, USA}
\begin{document}
%\ninept
%
\maketitle

\begin{abstract}
We {study} the problem of {constructing} a graph Fourier transform (GFT) for directed graphs (digraphs){, which decomposes graph signals into different modes of variation with respect to the underlying network. Accordingly, to capture low, medium and high frequencies we seek a digraph (D)GFT such that the orthonormal frequency components are as spread as possible in the graph spectral domain. This specification gives rise to a challenging nonconvex optimization problem, so we resort to a simple yet efficient heuristic to construct the DGFT basis from Laplacian eigenvectors of an undirected version of the digraph.} To {select frequency components which are as spread as possible}, we define a spectral dispersion function and show that it is supermodular. {Moreover, we show that orthonormality can be enforced via a matroid basis constraint, which motivates adopting a scalable greedy algorithm to obtain} an approximate solution with {provable} performance guarantee. The effectiveness of the {novel DGFT} is illustrated through numerical {tests} on synthetic and real-world graphs.
\end{abstract}
%
%%%%%%%%%%%%%%%%%%%%%%%%%%%%%%%%%%%%%%%%%%%
\section{Introduction}\label{S:Introduction}
\vspace{-2mm}
Network processes supported on the vertices of a graph can be viewed as graph signals, such as neural activities at different regions of the brain \cite{honey2007network}, or, infectious states of individuals in a population affected by an epidemic \cite{kolaczyk2014statistical}. Under the assumption that the signal properties relate to the underlying graph, the goal of graph signal processing (GSP) is to develop algorithms that fruitfully exploit this relational structure \cite{EmergingFieldGSP,sandryhaila2013}. {From this vantage point, } generalizations of traditional signal processing tasks such as filtering \cite{sandryhaila2013,teke2017extending,huang2016graph}, sampling and reconstruction \cite{marques2016sampling,chen2015discrete}, spectrum estimation \cite{marques2016stationary}, {(blind)} filter identification \cite{rasoul,segarra2017blindid} as well as signal representations \cite{thanou2014learning,zhu2012approximating} have been recently explored in the GSP literature.

An instrumental GSP tool is the graph Fourier transform (GFT), which decomposes a graph signal into orthonormal components describing different modes of variation with respect to the graph topology. Here we aim to generalize the GFT to directed graphs (digraphs); see also \cite{DSP_freq_analysis}, \cite{sardellitti_ICASSP}. We first propose a novel notion of signal variation (frequency) over digraphs and find an approximation of the maximum possible frequency ($f_{\text{max}}$) that a unit-norm graph signal can achieve. We design a digraph {(D)GFT} such that the resulting frequencies (i.e., the directed variation of the sought orthonormal bases) distribute as evenly as possible across $[0,f_{\text{max}}]$. To be more specific, we introduce some notations and basic GSP notions.

\noindent\textbf{Notation.} We consider a weighted digraph $\ccalG=(\ccalV,\bbA)$, where $\ccalV$ is the set of nodes {(i.e., vertices)} with cardinality $\lvert \ccalV \rvert=N$, and $\bbA \in \reals^{N \times N}$ is the graph adjacency matrix with entry $A_{ji}$ denoting the edge weight from node $i$ to  node $j$. We assume that the graph is connected and has no self loops; i.e. $A_{ii}=0$, and the edge weights are non-negative ($A_{ij} \geq 0$). We construct the underlying undirected graph $\ccalG^u=(\ccalV,\bbA^u)$ by replacing all directed edges of $\ccalG$ with undirected ones.
Let $\bbA^u \in \reals^{N \times N}$ denote the symmetric adjacency matrix of $\ccalG^u$, where by convention we set $A^u_{ij} = {A^u_{ji}}= \text{max}(A_{ij}, A_{ji})$. Then, the positive semi-definite Laplacian matrix takes the form $\bbL \coloneqq \bbD - \bbA^u$, where $\bbD$ is the diagonal degree matrix with $D_{ii} = \sum_{j} A^u_{ij}$. A graph signal $\bbx : \ccalV \mapsto \reals^N$ can be represented as a vector of size $N$, where component $x_i$ {denotes} the signal value at node $i \in \ccalV$.

\noindent\textbf{Related work.} For undirected graphs, the GFT of signal $\bbx$ {can be} defined as $\tilde{\bbx} = \bbV^T \bbx$, where $\bbV \coloneqq [\bbv_1, \ldots , \bbv_N]$ comprises the eigenvectors of the Laplacian \cite{EmergingFieldGSP}. Defining the total variation of the signal $\bbx$ with respect to the Laplacian $\bbL$ as
\begin{equation} \label{e:TV_def}
\text{TV}(\bbx) = \bbx^T \bbL \bbx = \sum_{i,j=1,j>i}^{N} A^u_{ij} (x_i - x_j)^2
\end{equation}
then it follows that the total variation of eigenvector $\bbv_k$ is $\text{TV}(\bbv_k) = \lambda_k$, the $k$\textsuperscript{th} Laplacian eigenvalue. Hence, eigenvalues {$0=\lambda_1<\lambda_2\leq\ldots\leq \lambda_N$} can be viewed as graph frequencies, indicating how the GFT bases vary over the graph. Note that there may be more than one eigenvector corresponding to a graph frequency in case of having repeated eigenvalues.
A more general GFT definition is based on the Jordan decomposition of adjacency matrix $\bbA=\bbV \bbJ \bbV^{-1}${,} where the frequency representation of graph signal $\bbx$ is $\tilde{\bbx} = \bbV^{-1} \bbx$ \cite{DSP_freq_analysis}. While valid for digraphs, the associated notion of signal variation in \cite{DSP_freq_analysis} does not ensure that constant signals have zero variation. Moreover, $\bbV$ is not necessarily orthonormal and Parseval's identity does not hold. {From a computational standpoint}, obtaining the Jordan decomposition is expensive and {often} numerically unstable; see also \cite{deri2017spectral}. Recently, a fresh look to the GFT for digraphs was put forth in \cite{sardellitti_ICASSP} based on minimization of {the} (convex) \lovasz extension of the graph cut size, subject to orthonormality constraints on the sought bases. However, the solution to such an optimization problem may not be unique. Also, the definition of cut size (and its \lovasz extension {which can be interpreted as a graph directed variation measure}) is based on {a} bipartition of the graph, while the network may have multiple (more than two) clusters. While the GFT bases in \cite{sardellitti_ICASSP} tend to be constant across clusters of the graph, in general they may fail to yield signal representations capturing different modes of signal variation with respect to $\ccalG$; see also Remark \ref{remark:star_graph}.

\noindent\textbf{Contributions.} Here instead we introduce a novel DGFT {(Section \ref{S:Notation})} that has the following desirable properties: P1) The bases offer notions of frequency and signal variation {overs} digraphs {which are} also consistent with {those used for} subsumed undirected graphs. P2) Frequencies are designed to be (approximately) equidistributed in $[0,f_{\text{max}}]$, to better capture low, middle, and high frequencies. P3) Bases are orthonormal so Parseval's identity holds and inner products are preserved in the vertex and graph frequency domains. Moreover, the inverse DGFT can be easily computed. To {construct a DGFT basis with the aforementioned} properties, a greedy algorithm is {outlined in Section \ref{S:fastdirected}} which is simple {(thus scalable to large graphs)} and efficient, while it offers provable performance guarantees.
%%%%%%%%%%%%%%%%%%%%%%%%%%%%%%%%%%%%%%%%%%%
\section{Preliminaries and Problem Statement }\label{S:Notation}
In this section we extend the notion of signal variation to digraphs and accordingly define graph frequencies. We then state the problem as {one} of finding an orthonormal basis with prescribed (approximately) equidistributed frequencies in the graph spectral domain.

Our goal is to find $N$ orthonormal bases capturing {different} modes of variation over the graph $\ccalG$. We collect these desired bases in a matrix $\bbU \coloneqq [\bbu_1, \ldots, \bbu_N] \in \reals^{N \times N}$, where $\bbu_i \in \reals^N$ represents the $i$\textsuperscript{th} frequency component. For undirected graphs, the quantity $\text{TV}(\bbx)$ in (\ref{e:TV_def}) measures how signal $\bbx$ varies over the {network} with Laplacian $\bbL$. This motivates defining a more general notion of signal variation for digraphs as {(cf. \cite[eq. (2)]{sardellitti_ICASSP})}
\begin{equation} \label{e:DV_def}
\text{DV}(\bbx) \coloneqq \sum_{i,j=1}^{N} A_{ji} [\bbx_i - \bbx_j]_{+}^2,
\end{equation}
where $[x]_{+} \coloneqq \text{max}(0,x)$. To { better appreciate \eqref{e:DV_def}}, consider a graph signal $\bbx \in \reals^N$ on digraph $\ccalG$ {and suppose a} directed edge represents the direction {of signal flow from a larger value to a smaller one. Thus,} an edge from node $i$ to node $j$ {(i.e., $A_{ji}>0$) contributes to $\text{DV}(\bbx)$ only if $x_i > x_j$. Notice that if $\ccalG$ is undirected, then $\text{DV}(\bbx)\equiv\text{TV}(\bbx)$.} 
Analogously to the undirected case, we define the frequency $f_i := \text{DV}(\bbu_i)$ as the directed variation of the basis $\bbu_i$.

{\noindent{\bf Problem statement.} Similar to the discrete spectrum of periodic time-varying signals, ideally we would like to have} $N$ equidistributed {graph} frequencies forming an arithmetic sequence
\begin{equation} \label{e:directed_freq}
f_k = {\text{DV}(\bbu_k)=} \frac{k-1}{N-1} f_{\text{max}}, \quad k=1,\ldots,N
\end{equation}
where $f_\text{max}$ is the maximum variation of a unit-norm graph signal on $\ccalG$. {Accordingly, given $\ccalG$ and for prescribed graph spectrum $\{f_k\}_{k=1}^N$, the DGFT basis $\bbU$ can be found by solving the following non-convex optimization problem ($\bbI_N$ is the $N\times N$ identity matrix)}
\begin{equation} \label{e:opt_prob}
\bbU= \argmin_{\bbU \in \reals^{N\times N}}
\sum_{k=1}^N (\text{DV}(\bbu_k) - f_k)^2,\quad \text{s. t. }
\bbU^T\bbU=\bbI_N.
\tag{\ccalP1}
\end{equation}
{Problem \eqref{e:opt_prob} can be tackled by splitting methods for orthogonality constrained problems, e.g. \cite{lai2014splitting}. However those methods are computationally expensive and do not offer convergence guarantees for non-convex objectives such as in \eqref{e:opt_prob}. In Section \ref{S:fastdirected}, we will propose a simple yet efficient heuristic to construct the DGFT basis from Laplacian eigenvectors of $\ccalG^u$.}

{Going back to design of the graph spectrum $\{f_k\}_{k=1}^N$, a few noteworthy challenges remain.} {First}, attaining the exact frequencies in \eqref{e:directed_freq} may be impossible. {This can be clearly seen for undirected graphs, where one has the additional constraint that the summation of frequencies is constant, since}  $\sum_{k=1}^{N}f_k=\sum_{k=1}^{N}\text{TV}(\bbu_k)=\text{trace}(\bbL)$.

{Second}, one needs to determine the maximum frequency $f_{\max}$ that a unit-norm basis can attain. In $\ccalG^u$, it is immediate that
\begin{equation} \label{e:lambda_max}
f_{\text{max}}^{u} = \max_{\lVert \bbu \rVert = 1} \text{TV}(\bbu) =\max_{\lVert \bbu \rVert = 1} \bbu^T \bbL \bbu = \lambda_{\text{max}},
\end{equation}
where $\lambda_{\text{max}}$ is the largest eigenvalue of the Laplacian matrix $\bbL$.
However, finding the maximum directed variation is in general challenging, since one needs to solve {the (non-convex) spherically-constrained problem}
\begin{equation}\label{e:opt_f_max}
\bbu^\star = \argmax_{\lVert \bbu \rVert = 1} \quad \text{DV}(\bbu) \quad \text{and} \quad f_{\text{max}} := \text{DV}(\bbu^\star).
\end{equation}
%
%Unfortunately, there is no efficient way to maximize a general convex function over a spherical constraint. 
Still, it follows that $f_{\text{max}}$ can be upper-bounded by $\lambda_{\text{max}}$. This is because dropping the direction of any edge can not decrease the directed variation, hence $\lambda_{\text{max}}$ gives an upper bound for the maximum DV. In the next section, we propose a way to find a basis with variation of at least $f_{\text{max}}/2$. But before moving on, a remark is in order.
%A heuristic method for solving the spherical constrained problems based on Bregman iteration is proposed in \cite{lai2014splitting}. 

\begin{figure}[t]
	\centering    
	{\includegraphics[width=0.7\linewidth]{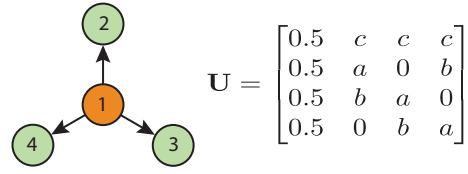}}
	%\subfigure[multiprimary]{\includegraphics[width=0.30\linewidth]{\FigDirPresenta/Multiprimary_example}}
	%\subfigure[$CIELUV$]{\includegraphics[width=0.45\linewidth]{\FigDir/REC709_GamutLUV}}        
	\caption{A toy digraph and {the GFT basis $\bbU$ from~\cite{sardellitti_ICASSP}, where $a=(\sqrt{5}+1)/4$, $b=(1-\sqrt{5})/4$ and $c=-0.5$}.}
	\label{fig:example}
\end{figure}

\begin{remark}[{Motivation for} spread frequencies] \label{remark:star_graph}\normalfont
Consider {the} toy graph with $4$ nodes shown in Fig. \ref{fig:example} (left).
Let $\mathrm{DV}'(\bbx) \coloneqq \sum_{i,j=1}^{N} A_{ji} [\bbx_i - \bbx_j]_{+}$ be a directed variation measure, {and consider minimizing $\sum_{k=1}^N \mathrm{DV}'(\bbu_k),\:\text{ s. t. }
	\bbU^T\bbU=\bbI_N$ as in \cite{sardellitti_ICASSP}}.  {The optimum basis $\bbU$ is shown} in Fig. \ref{fig:example} (right),
%\begin{equation}
%\bbU=
%\begin{bmatrix}
%0.5 & c  & c & c\\
%0.5 & a  & 0 & b\\
%0.5 & b  & a & 0\\
%0.5 & 0  & b & a
%\end{bmatrix},
%\end{equation}
where $a=(\sqrt{5}+1)/4\approx 0.8090$, $b=(1-\sqrt{5})/4\approx -0.3090$, and $c=-0.5$. These values are chosen such that $a+b+c=0$, $a^2+b^2+c^2=1$, and $c^2+ab=0$ which implies the orthonormality of $\bbU$. {As a result,} all columns of $\bbU$ satisfy $\mathrm{DV}'(\bbu_k)=0,\: k=1,\ldots, 4,$ and the {synthesis} formula $\bbx = \bbU \tilde{\bbx}$ {fails to} offer an expansion of {\bbx} with respect to \emph{different} modes of variation {(e.g., low and high graph frequencies)}.
\end{remark}
%%%%%%%%%%%%%%%%%%%%%%%%%%%%%%%%%%%%%%%%%%%
\section{{A Digraph Fourier Transform Heuristic}}\label{S:fastdirected}
%\begin{mydefinition}[Matroid]
%Let $S$ be a finite ground set, and let $\ccalI$ be a collection of subsets of $S$. The pair $\ccalM=(S,\ccalI)$ is a matroid if the following conditions hold: 
%\begin{itemize}
%\item Hereditary Property: If $T\in \ccalI$, then $T'\in \mathcal{I}$ for all $T'\subseteq T$.
%\item Augmentation Property: If $T_1,T_2\in \mathcal{I}$ and $\abs{T_1}<\abs{T_2}$, then there exists $e\in T_2\backslash T_1$ such that $T_1\cup \{e\}\in \mathcal{I}$.
%\end{itemize}
%The collection $\mathcal{I}$ is called the set of independent sets of the matroid $\mathcal{M}$. A maximal independent set is a basis. It is easy to show that all the bases of a matroid have the same cardinality.
%\end{mydefinition}
%Note that the bases of the partition matroid are those sets $A$ that satisfy $\abs{A\cap S_i}=d_i$, for $i=1,...,m$.

As mentioned in {Section} \ref{S:Notation}, one challenge in finding a well distributed set of frequencies on a digraph is to calculate the maximum frequency $f_{\text{max}}$ [cf. (\ref{e:opt_f_max})]. Here we show how to find a basis vector $\tilde{\bbu}$ with an approximate $\tilde{f}_{\text{max}}\coloneqq \text{DV}(\tilde{\bbu})$ which is at least half of $f_{\text{max}}$.
\begin{myproposition}
For a digraph $\ccalG$, {recall $\ccalG^{u}$ and the spectral radius $\lambda_{\text{max}}$ of its Laplacian $\bbL$. Let $\bbu$ be the dominant eigenvector.} Then,
\begin{equation}
\tilde{f}_{\mathrm{max}}\coloneqq \max{\{\mathrm{DV}(\bbu),\mathrm{DV}(-\bbu)\}}\geq \frac{f_{\mathrm{max}}}{2}.
\end{equation} 
\end{myproposition}
\begin{myproof}
	First recall that
	\begin{equation} \nonumber
	\lambda_{\text{max}}=\bbu^T\bbL\bbu= \sum_{i,j=1,j>i}^{N} A^u_{ij} (u_i - u_j)^2=\frac{1}{2}\sum_{i,j=1}^{N} A^u_{ij} (u_i - u_j)^2.
	\end{equation}
	Since $A^u_{ij}\leq A_{ij}+A_{ji}$ and $(u_i - u_j)^2=[u_i - u_j]_{+}^2+[u_j - u_i]_{+}^2$, then
	\begin{align} \nonumber
	\lambda_{\text{max}}\leq \frac{1}{2}\sum_{i,j=1}^{N} (A_{ij}+A_{ji}) \bigg([u_i - u_j]_{+}^2+[u_j - u_i]_{+}^2\bigg)\nonumber \\
	=\frac{1}{2} \big(2\times \text{DV}(\bbu)+2\times \text{DV}(-\bbu)\big)=\text{DV}(\bbu)+\text{DV}(-\bbu).\nonumber
	\end{align}
{In conclusion}, at least one of $\text{DV}(\bbu)$ or $\text{DV}(-\bbu)$ is larger than $\lambda_{\text{max}}/2$, and this completes the proof since $\lambda_{\text{max}}\geq f_{\text{max}}$.
\end{myproof}

\noindent In practice, we can {compute} $\max{\{\text{DV}(\bbu),\text{DV}(-\bbu)\}}$  {for} \emph{any}  eigenvector {$\bbu$} of the Laplacian matrix. This will possibly give a higher frequency {in $\ccalG$}, while preserving the $1/2$-approximation.

Fixing $f_1=0$ and $f_N=\tilde{f}_{\text{max}}$, {in lieu of \eqref{e:directed_freq} we will henceforth construct} a disperse set of frequencies by using the eigenvectors of $\bbL$. Let $f_i\coloneqq \text{DV}(\bbu_i)$ and $\overline{f}_i\coloneqq \text{DV}(-\bbu_i)$, where $\bbu_i$ is the $i$\textsuperscript{th} eigenvector of $\bbL$. Define the set of all {candidate} frequencies {as} $F:=\{f_i,\overline{f}_i:1<i<N\}$. {The goal is to select} $N-2$ frequencies from $F$ such that together with $\{0,\tilde{f}_{\text{max}}\}$ they form our Fourier frequencies. To preserve orthonormality, we would opt exactly one from each pair $\{f_i,\overline{f}_i\}$. We will argue later that this induces a matroid basis constraint.\vspace{0.1cm}

{\noindent{\bf Optimization scheme.}} {To find the DGFT basis, we define a spectral dispersion function} that measures how well {spread are} the {corresponding} frequencies {over} $[0,\tilde{f}_{\text{max}}]$. For {frequency set} $S\subseteq F$, let $s_1\leq s_2\leq ...\leq s_m$ be the elements of $S$ in non-decreasing order, where $m=\abs{S}$. Then we define the dispersion of $S$ as
\begin{equation} \label{e:dispersion_def}
\delta(S)=\sum_{i=0}^{m} (s_{i+1}-s_i)^2,
\end{equation}
where $s_0=0$ and $s_{m+1}=\tilde{f}_{\text{max}}$.
%\begin{equation}
%\delta(S)=\sup_{x\in [0,\tilde{f}_{\text{max}}]} \min_{f\in S\cup\{0,\tilde{f}_{\text{max}}\}} \abs{x-f}
%\end{equation}
It can be verified that $\delta(S)$ is a monotone non-increasing function, which means that for any sets $S_1\subseteq S_2$, we have $\delta(S_1)\geq \delta(S_2)$. For a fixed value of $m$, one can show that $\delta(S)$ is minimized when the $s_i$'s form an arithmetic sequence, {consistent} with our design goal. Hence, we {seek} to minimize $\delta(S)$ through a set function optimization procedure. In Lemma \ref{lemma:obj}, we show that the dispersion function (\ref{e:dispersion_def}) has the supermodular property. First, we define submodularity/supermodularity.
\begin{mydefinition}[Submodularity]
Let $S$ be a finite ground set. A set function $f:2^S\mapsto \reals$ is submodular if:
\begin{equation}
\label{eq:def_submod}
f(T_1\cup \{e\})-f(T_1)\geq f(T_2\cup \{e\})-f(T_2),
\end{equation}
for all subsets $T_1\subseteq T_2\subseteq S$ and any element $e\in S\backslash T_2$.
\end{mydefinition}
A set function $f$ is said to be supermodular if $-f$ is submodular, i.e. \eqref{eq:def_submod} holds in {the} other direction. Submodularity captures the \emph{diminishing returns} property that arises in many applications such as facility location, sensor placement, feature selection \cite{nemhauser1978analysis}.
%%%%%%%%%%%%%%%%%% Lemma 1 %%%%%%%%%%%%%%%%%%
\begin{mylemma}
\label{lemma:obj}
The {spectral} dispersion function $\delta:2^F\mapsto \reals$ defined in (\ref{e:dispersion_def}) is a supermodular function.
%i.e. adding a single element $e$ to a superset $S_2$, produces bigger marginal value compared to adding the same element to a subset $S_1\subseteq S_2$. More precisely, we have to show that the following equation holds for any sets $S_1\subseteq S_2\subseteq F$ and any element $e\in F\backslash S_2$.
%\begin{equation}
%\delta(S_1\cup \{e\})-\delta(S_1)\leq \delta(S_2\cup \{e\})-\delta(S_2). \label{eq:dispersionsupermodularity}
%\end{equation}
\end{mylemma}
\begin{myproof}
Consider two subsets $S_1,S_2$ such that $S_1\subseteq S_2\subseteq F$, and a single element $e\in F\backslash S_2$. Let $s_1^L$ and $s_1^R$ be the largest value smaller than $e$ and smallest value greater than $e$ in $S_1\cup\{0,\tilde{f}_{\text{max}}\}$, respectively {(i.e.,} $e$ breaks the gap between $s_1^L$ and $s_1^R$). Similarly, let $s_2^L$ and $s_2^R$ be defined for $S_2$. Since $S_1\subseteq S_2$, then $s_1^L\leq s_2^L\leq e\leq s_2^R\leq s_1^R$. {The result follows by comparing} the marginal values
	\begin{align*}
	&\delta(S_1\cup \{e\})-\delta(S_1)=(s_1^R-e)^2+(e-s_1^L)^2-(s_1^R-s_1^L)^2\\
	&=-2(s_1^R-e)(e-s_1^L)\leq -2(s_2^R-e)(e-s_2^L)\\
	%&=(s_2^R-e)^2+(e-s_2^L)^2-(s_2^R-s_2^L)^2\\
	&=\delta(S_2\cup \{e\})-\delta(S_2).
	\end{align*}
\end{myproof}

%\vspace{-0.5cm}
%%%%%%%%%%%%%%%%%%%%%%%%%%%%%%
\begin{figure*}[t]
	\centering
	\def \thisplotscale {0.95}
\def \unit {\thisplotscale cm}

{\small \begin{tikzpicture}[thick, x = 1*\unit, y = 1*\unit]

\node (graph_fig_1_a) at (0,-2.5)
{\includegraphics[width=0.292\textwidth]{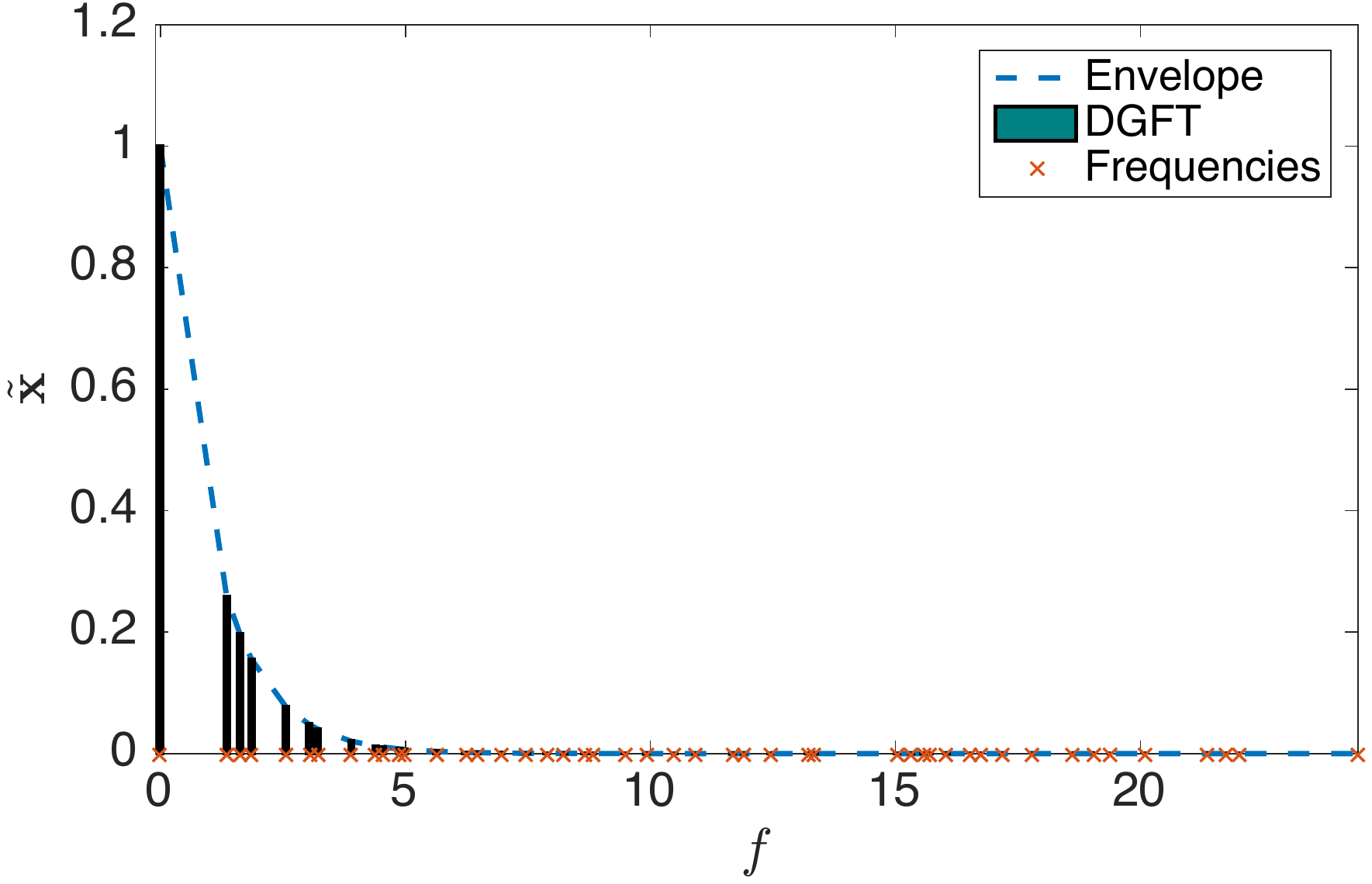}}; % HERE PUT FIG 1 A

%\node (graph_fig_1_b) at (0,-3.5)
%{\includegraphics[width=0.30\textwidth, height=0.1\textwidth]{figures/noisy_log.pdf}}; % HERE PUT FIG 1 B

\node (a) at (0,-4.6) {{{\small (a)}}};

\node (graph_fig_2) at (6,-2.5)
{\includegraphics[width=0.303\textwidth]{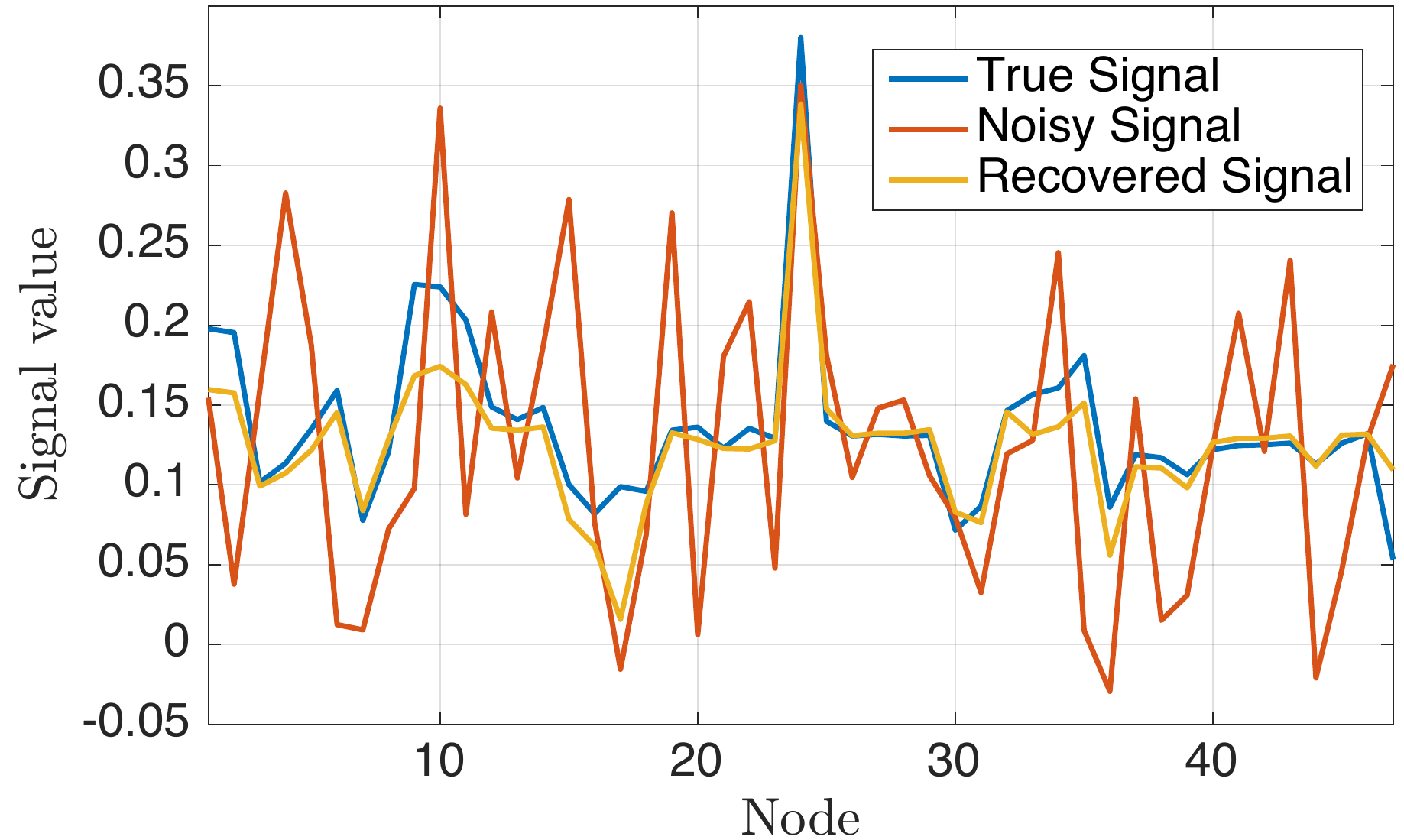}}; % THIS IS THE FIGURE FOR SIM 2

\node (b) at (6,-4.6) {{{\small (b)}}};

\node (graph_phase_3) at (12,-2.47)
{\includegraphics[width=0.295\textwidth]{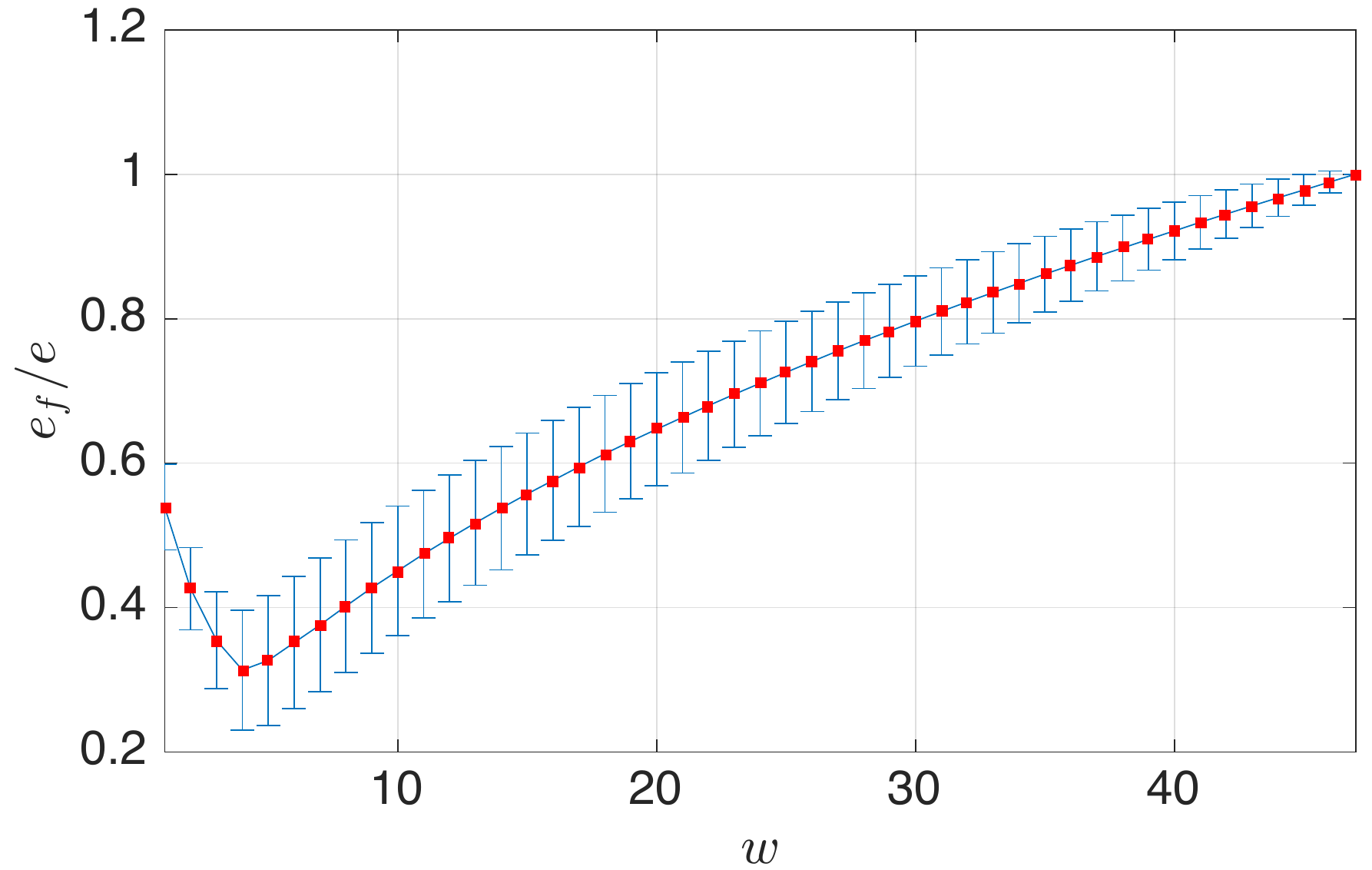}}; % HERE PUT FIG 3

\node (c) at (12,-4.6) {{{\small (c)}}};

\end{tikzpicture}}
	\vspace{-0.2cm}
	\caption{(a) A synthetic signal defined in the graph spectral domain by ${\tilde{x}_f} =ce^{-f}$. (b) True, noisy, and recovered signal values at different nodes for a sample realization. (c) Average recovery error using windowed filter over recovery error without filtering versus the window size.}
	\vspace{-0.15in}
	\label{F:num_exp}
\end{figure*}
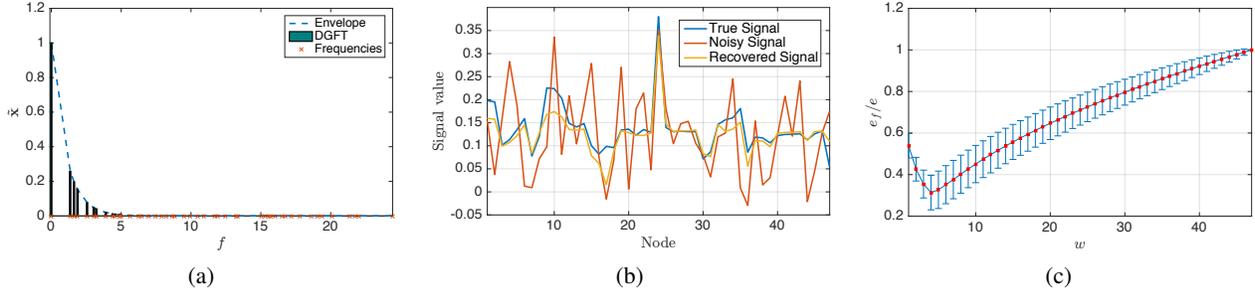
%%%%%%%%%%%%%%%%%%%%%%%%%%%%%%
\noindent {Recalling the orthonormality constraint, we define} $\ccalB$ to be the set of all subsets $S\subseteq F$ that satisfy $\abs{S\cap \{f_i,\overline{f_i}\}}=1$, $i=2,...,N-1$. Then, {frequency selection from $F$ boils down to solving}
\begin{equation}
 {\min_{S}  \delta(S), \quad 
 \text{s. t. } S\in \ccalB.}
\tag{\ccalP2}\label{eq:supmin}
\end{equation}
{Next}, in Lemma \ref{lemma:constraint} we show that the constraint in \eqref{eq:supmin} is a matroid basis constraint. {But first, we define the notion of partition matroid.}
\begin{mydefinition}[Partition matroid \cite{schrijver}]
\label{def:partition_matroid}
Let $S$ denote a finite set, and let $S_1,...,S_m$ be a partition of $S$, i.e. a collection of disjoint sets such that $S_1\cup ... \cup S_m=S$. Let $d_1,...,d_m$ be a collection of non-negative integers. Define a set $\ccalI$ by $A\in \ccalI$ iff $\abs{A\cap S_i}\leq d_i$ for all $i=1,...,m$. Then, $\ccalM=(S,\ccalI)$ is called the partition matroid.
\end{mydefinition}
%
%%%%%%%%%%%%%%%%%% Lemma 2 %%%%%%%%%%%%%%%%%%
\begin{mylemma}
\label{lemma:constraint}
There exists a (partition) matroid $\mathcal{M}$ such that the set $\mathcal{B}$ in (\ref{eq:supmin}) is the set of all bases of $\mathcal{M}$.
\end{mylemma}
\begin{myproof}
{Recall Definition \ref{def:partition_matroid} and set $S=F$, $S_i=\{f_i,\overline{f_i}\}$ and $d_i=1$ for all $i=2,...,N-1$}, to get a partition matroid $\ccalM=(F,\ccalI)$. Moreover, the bases of $\ccalM$, which are defined as the maximal elements of $\ccalI$, are those subsets $A\subseteq F$ that satisfy $\abs{A\cap \{f_i,\overline{f_i}\}}=1$ for all $i=2,...,N-1$, which are the elements of $\ccalB$.
\end{myproof}

\noindent{Lemmas} \ref{lemma:obj} and \ref{lemma:constraint} {assert} that \eqref{eq:supmin} is a supermodular minimization problem subject to a matroid basis constraint. Since supermodular minimization is NP-hard and hard to approximate to any factor \cite{kelner2007hardness,mittal2013fptas}, we create a submodular function $\tilde{\delta}(S)$ and use the algorithms for submodular maximization to find a set of disperse bases ${\bbU}$. In particular, we define
\begin{equation}
\tilde{\delta}(S):=\tilde{f}_{\text{max}}^2 - \delta(S),
\end{equation}
{which} is a non-negative (increasing) submodular function, because $\delta(\emptyset)=\tilde{f}_{\text{max}}^2$ is an upper bound for $\delta(S)$. There are several results for maximizing submodular functions under matroid constraints for {both the} non-monotone \cite{lee2009non} and monotone cases \cite{fisher1978analysis,calinescu2011maximizing}. We {adopt} the greedy algorithm of \cite{fisher1978analysis} due to simplicity ({tabulated} under Algorithm \ref{alg:greedy}), which offers the following performance guarantee.
%=================================== ALGORITHM
%\renewcommand\algorithmicdo{}	% removes "DO" from for loops
\begin{algorithm}[t]
\caption{Greedy Dispersion Minimization}
\label{alg:greedy}
\begin{algorithmic}[1]
	\STATE \textbf{Input:} Set of possible frequencies $F$.
	\STATE \textbf{Initialize} $S=\emptyset$.
	\REPEAT
		\STATE $e \gets \argmax_{f\in F} \left\{\tilde{\delta}(S\cup\{f\})-\tilde{\delta}(S)\right\}$.
		\STATE $S \gets S\cup \{e\}$.
		\STATE Delete from $F$ the pair $\{f_i,\overline{f}_i\}$ that $e$ belongs to.
	\UNTIL{$F=\emptyset$}
\end{algorithmic}
\end{algorithm}
%=================================	
\begin{mytheorem} [\cite{fisher1978analysis}] \label{theorem:approximation}
 Let $S^*$ be the solution of problem (\ref{eq:supmin}) and $S^{\text{g}}$ be the output of the greedy Algorithm \ref{alg:greedy}. Then,
\begin{equation}
\tilde{\delta}(S^{\text{g}})\geq\frac{1}{2}\times\tilde{\delta}(S^*).
\end{equation}
\end{mytheorem}
Notice that Theorem \ref{theorem:approximation} offers a worst-case guarantee, and Algorithm \ref{alg:greedy} is usually able to find near-optimal solutions in practice.\vspace{0.1cm}

All in all, the DGFT {basis construction} algorithm entails the following steps. First, {we form $\ccalG^u$ and find} the eigenvectors of $\bbL$. Second, the set $F$ is formed by calculating DV for each eigenvector {$\bbu_i$} and its negative {$-\bbu_i$, $i=2,\ldots,N-1$.} Finally, the greedy Algorithm \ref{alg:greedy} is run on the set $F$, and the output determines the set of frequencies as well as the orthonormal set of DGFT bases.
	
%%%%%%%%%%%%%%%%%%%%%%%%%%%%%%%%%%%%%%%%%%%
\section{Numerical Results}\label{S:numerical}

\begin{figure}[t]
	\centering    
	{\includegraphics[width=0.95\linewidth]{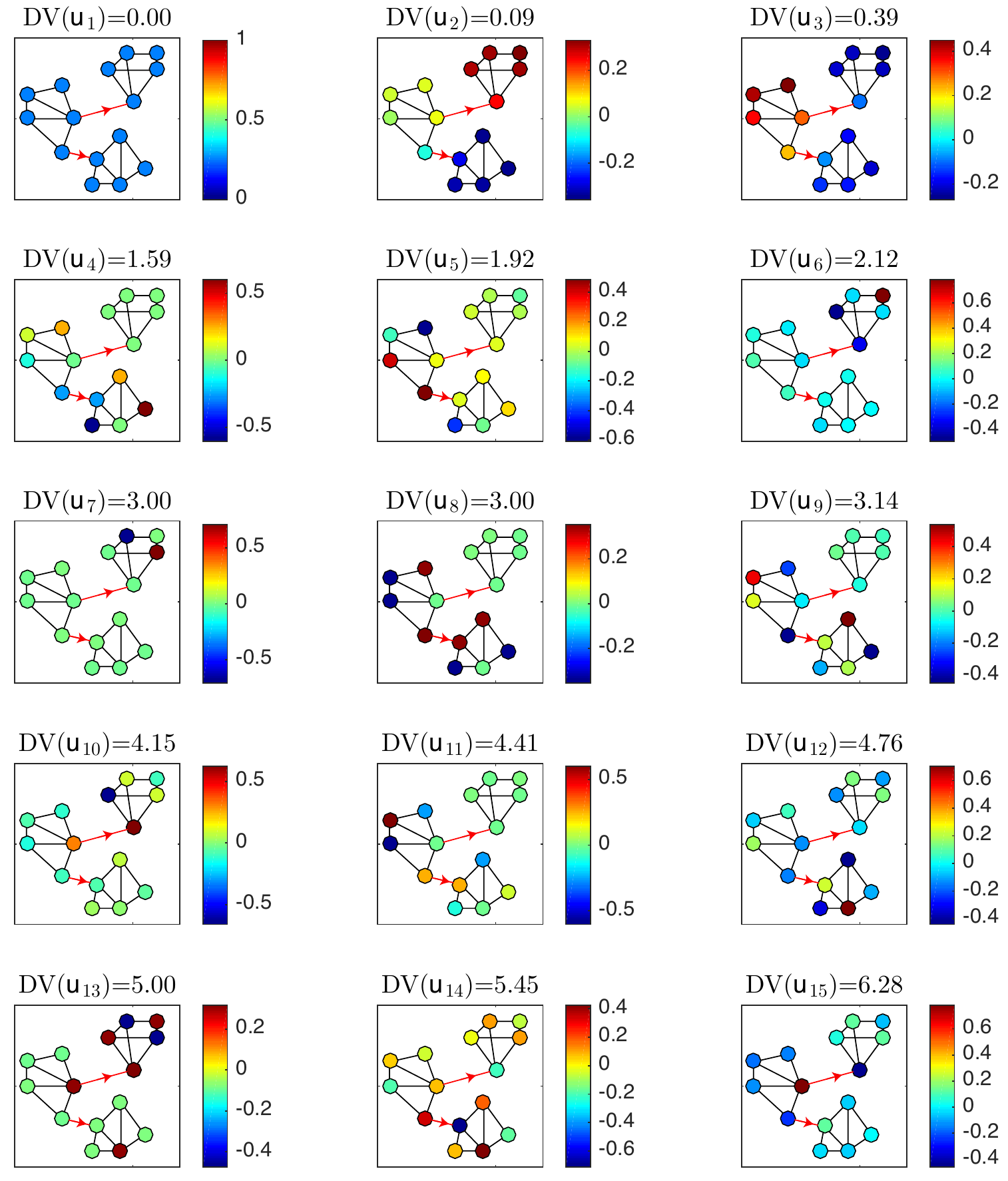}}
	%\subfigure[multiprimary]{\includegraphics[width=0.30\linewidth]{\FigDirPresenta/Multiprimary_example}}
	%\subfigure[$CIELUV$]{\includegraphics[width=0.45\linewidth]{\FigDir/REC709_GamutLUV}}   
%	\vspace{-0.3cm}     
	\caption{DGFT bases obtained using Algorithm \ref{alg:greedy}, along with their respective DV values (frequencies).}
	\label{fig:num_1}
\end{figure}
Here we study the performance of Algorithm \ref{alg:greedy} to construct the DGFT basis, via simulations on two graphs.

\noindent\textbf{Synthetic graph.} First, we construct the DGFT for a digraph with $N = 15$ nodes and compare it with the GFT in \cite{sardellitti_ICASSP}. Fig. \ref{fig:num_1} shows all the bases and their frequencies derived from Algorithm \ref{alg:greedy} and (\ref{e:DV_def}), respectively. The DGFT can better capture low, middle, and high frequencies in comparison with the variations exhibited by the bases learned via algorithm in \cite{sardellitti_ICASSP}; see Table \ref{tab:1}. For fairness of comparison, we calculate the directed variation of the proposed bases using $\text{DV}'${; see Remark \ref{remark:star_graph}.} Apparently, the DGFT bases {span} a wider range of {variations}. Moreover, we rescale DV$'$ values in Table \ref{tab:1} to the $[0,1]$ interval and calculate their dispersion using (\ref{e:dispersion_def}). The dispersion of the DGFT bases is $0.16$ and the counterpart in \cite{sardellitti_ICASSP} gives a dispersion of $0.24$, which confirms the proposed method yields a better (i.e. smaller) frequency spread.

%\begin{table}[t]
%
%	\centering
%		\footnotesize
%	\begin{tabular}{@{} |>{\centering}m{0.35cm}|c|c|>{\centering}m{0.35cm}|c|c|>{\centering}m{0.35cm}|c|c| @{}}
%		\hline
%		Basis & DV$'$ & DV$'$ & Basis & DV$'$ & DV$'$ & Basis & DV$'$ & DV$'$\\
%		     & Alg. \ref{alg:greedy} & \cite{sardellitti_ICASSP} & & Alg. \ref{alg:greedy} & \cite{sardellitti_ICASSP}  & & Alg. \ref{alg:greedy} & \cite{sardellitti_ICASSP} \\
%		\hline
%		1	&      0      &	0	  &	6   &		4.35 &   	3	    	&	11  	& 6.55 & 4.08	 \\	
%		2	&	0.97 &    0		  &    7   &	  	4.90 &   	3.37		&	12  	& 7.78 & 4.25	  \\    
%		3	&	2.19 &    0		  &    8   &	  	5.66 &   	3.46		&	13	& 7.90 & 4.62    \\
%		4	&	4.12 & 2.24	  &	9   &		6.01 &   	3.53 		&	  14  	& 8.09 & 4.62	\\
%		5	&	4.25 & 2.46	  &	10 &   	6.53 &	3.67 		&	 15  	& 8.70 & 4.98	  \\
%		\hline
%	\end{tabular}
%	\caption{directed variations (frequencies) of obtained bases.}
%	\label{tab:1}
%\end{table}

\begin{table}[t]
	
	\centering
	\footnotesize
	\begin{tabular}{|c|c|c||c|c|c|}
		\hline
		Basis & DV$'$ & DV$'$ & Basis & DV$'$ & DV$'$\\
		& Alg. \ref{alg:greedy} & \cite{sardellitti_ICASSP} & & Alg. \ref{alg:greedy} & \cite{sardellitti_ICASSP}  \\
		\hline
		1	&      0   &	0	   &	9   &		6.01 &   	3.53	    	 \\	
		2	&	0.97 &    0		 &    10 &   	6.53 &	3.67		  \\    
		3	&	2.19 &    0		 &   11  	& 6.55 & 4.08		   \\
		4	&	4.12 & 2.24	   &	12  	& 7.78 & 4.25 			\\
		5	&	4.25 & 2.46	   &	13	& 7.90 & 4.62 		 \\
		6   &   4.35 &   3    &	  14  	& 8.09 & 4.62 \\
		7   &	4.90 &   3.37	&	 15  	& 8.70 & 4.98 \\
		8   &	5.66 &   3.46 & &  &\\
		\hline
	\end{tabular}
	\caption{{Directed} variations (frequencies) of obtained bases.}
	\label{tab:1}
\end{table}

\noindent\textbf{Brain graph.} Next, we consider a real brain graph to demonstrate the effectiveness of {the} DGFT in a denoising task. This graph represents an anatomical connection network of the macaque cortex. It consists of 47 nodes and 505 edges (among which 121 links are directed), which is used e.g. in \cite{honey2007network,rubinov2010complex}.
%This network includes the ventral and dorsal streams of visual cortex, and groups of somatosensory and somatomotor regions.
The vertices represent different hubs in the brain, and the edges capture directed information flow among them.
Let $\bbU$ be the orthonormal DGFT basis obtained via Algorithm \ref{alg:greedy}. We construct a smooth graph signal in the spectral domain as ${\tilde{x}_f} = ce^{-f}$, where $c$ is calculated such that the graph signal $\bbx = \bbU \tilde{\bbx}$ has unit norm. The DGFT of signal $\bbx$ is shown in Fig. 2 (a). We aim to recover the signal from noisy measurements $\bby = \bbx {+} \bbn$, where the additive noise $\bbn$ is a zero-mean, Gaussian random vector with covariance matrix $10^{-2} \bbI_{N}$. To that end, we use a low-pass filter $\tilde{\bbH} \coloneqq \text{diag}(\tilde{\bbh})$, where $\tilde{h}_{i} = \mathbb{1} [i \leq w]$ and $w$ is the spectral window size. The filter selects the first $w$ components of the signal's DGFT, and we approximate the noisy signal by
\begin{equation} \label{e:filtered_signal}
\hat{\bbx} = \bbU \tilde{\bbH}  \tilde{\bby} = \bbU \tilde{\bbH} \bbU^T \bby.
\end{equation}
Fig. 2(b) shows a sample noisy signal ($\bby$) along with the recovered signal ($\hat{\bbx}$) using (\ref{e:filtered_signal}) with $w=4$, and the original signal ($\bbx$). Filter design and {choice of} $w$ is out of scope of this paper, but as we can see $\hat{\bbx}$ closely approximates $\bbx$. Moreover, we compute the relative recovery error with and without low-pass filtering as $e_f = \lVert \hat{\bbx} - \bbx \rVert/\lVert \bbx \rVert$ and $e=\lVert \bbn \rVert/\lVert \bbx \rVert$, respectively. Fig. 2(c) depicts $e_f/e$ versus $w$ averaged over $1000$ Monte-Carlo simulations, and demonstrates the effectiveness of adopting filters along with the proposed DGFT.
\vspace{-0.2cm}
%%%%%%%%%%%%%%%%%%%%%%%%%%%%%%%%%%%%%%%%%%
\section{Conclusion}\label{S:Conclusions}
In this paper, we introduced a new method to find an orthonormal set of graph Fourier bases for digraphs. To that end, we proposed a measure of directed variation to capture the notion of frequency. Our approach is to find a basis that spans the entire frequency range and for which frequency components are as evenly distributed as possible. We showed that the maximum directed variation can be approximately achieved, and then used the eigenvectors of the Laplacian matrix of the underlying undirected graph to find the Fourier basis. We defined a dispersion function to quantify the quality of any feasible solution compared to our ideal design, and proposed a fast greedy algorithm to minimize this dispersion. The greedy algorithm offers theoretical approximation guarantees by virtue of matroid theory and results for submodular function optimization.

{With regards to future directions, }the complexity of finding {the} maximum frequency ($f_\text{max}$) on a digraph is an interesting open question. If NP-hard, it will be interesting to find the best achievable approximation factor (here we gave a $1/2$-approximation). Furthermore, our performance guarantee measures the quality of the proposed basis with respect to the optimal set one can get by using the eigenvectors of the Laplacian matrix. A significant improvement {would be to} generalize this guarantee to any orthonormal basis, or {otherwise establish} that the global minimizer of dispersion (among all orthonormal bases) lies within the considered feasible set $\ccalB$. 

\vfill\pagebreak

%\clearpage
%\newpage
%\bibliographystyle{IEEEbib}
%\bibliography{refs}

\bibliographystyle{IEEEtran}
\bibliography{references}

\end{document}